\documentclass[runningheads,a4paper]{llncs}

\usepackage{graphicx}

\usepackage{url}
\urldef{\mail}\path|majdisz@kiv.zcu.cz|
\urldef{\mails}\path|smolik@kiv.zcu.cz|

\newcommand{\keywords}[1]{\par\addvspace\baselineskip
\noindent\keywordname\enspace\ignorespaces#1}

\usepackage[hidelinks]{hyperref}
\makeatletter \let\cl@chapter\relax \makeatother
\usepackage[capitalise,nameinlink]{cleveref}
\crefname{section}{Sect.}{Sections}
\Crefname{section}{Section}{Sections}
\crefname{figure}{Fig.}{Figures}
\Crefname{figure}{Figure}{Figures}
\crefname{equation}{}{}
\crefname{table}{Table}{Tables}

\usepackage[caption=false]{subfig}
\usepackage{bm}
\usepackage{mathtools}
\usepackage{amsmath, amsfonts, amssymb}
\usepackage{xfrac}
\usepackage{multirow}
\usepackage{csquotes}

\usepackage{algorithm}
\usepackage[noend]{algpseudocode}

\crefname{algorithm}{Algorithm}{Algorithms}

\usepackage{dblfloatfix}
\usepackage{arydshln}

\usepackage{xcolor} 
\usepackage{setspace}

\usepackage{microtype} 

\begin{document}
\mainmatter  
%
\title{Determination of Stationary Points and Their Bindings in Dataset using RBF Methods}
%
\titlerunning{Determination of Stationary Points and Their Bindings}
%
\author{Zuzana Majdisova%
\thanks{Corresponding author.}%
\and Vaclav Skala\and Michal Smolik}
\authorrunning{Determination of Stationary Points and Their Bindings}
%
\institute{Department of Computer Science and Engineering, Faculty of Applied Sciences, University of West Bohemia, \\
Univerzitn\'{i} 8, CZ 30614 Plze\v{n}, Czech Republic\\
\mail, \mails\\
\url{www.vaclavskala.eu}}
%
%
\toctitle{Determination of Stationary Points and Their Bindings in Dataset using RBF Methods}
\tocauthor{Z. Majdisova et al.}
\maketitle

\begin{abstract}
Stationary points of multivariable function which represents some surface have an important role in many application such as computer vision, chemical physics, etc. Nevertheless, the dataset describing the surface for which a sampling function is not known is often given. Therefore, it is necessary to propose an approach for finding the stationary points without knowledge of the sampling function.\\

In this paper, an algorithm for determining a set of stationary points of given sampled surface and detecting the bindings between these stationary points (such as stationary points lie on line segment, circle, etc.) is presented. Our approach is based on the piecewise RBF interpolation of the given dataset.
\end{abstract}

\section{Introduction}
Stationary points of the given explicit function $f(\bm{x})$ are points where the gradient of the function $f(\bm{x})$ is zero in all directions, i.e. all partial derivatives are zero:
\begin{equation}
\label{eq:stacPts}
\begin{aligned}
\displaystyle \nabla f(\bm{x}) &= \bm{0} \qquad \bm{x}\in \mathbb{E}^n \textrm{, i.e.}\\
\displaystyle \frac{\partial f(\bm{x})}{\partial x_k} &= 0 \qquad k = 1,\dots, n\textrm{,}
\end{aligned}
\end{equation}
\noindent where $n$ denotes the dimension of space. The knowledge of stationary points is required in many areas that are used a multidimensional data analysis, e.g. \cite{Banerjee1985}, \cite{Tsai1993}, \cite{Comaniciu02}, \cite{Strodel2008}, \cite{Liu2011}. The significant features of the given dataset can be determined using the set of stationary points. This properties can be further used for improving the quality of the RBF approximation \cite{Majdisova2017}, \cite{Majdisova2017b}, etc.
In the technical applications, the sampling function is not often known and only the dataset describing the given surface is specified. Therefore, it is necessary determining the stationary points without knowledge of the sampling function. Moreover, for a higher dimension of space $n\ge 2$, it is possible that the stationary points of given surface are not only isolated but they can be formed into line segments, circles or some other shapes. A new approach for searching of bindings between stationary points will be described in this paper. Knowledge of these bindings is suitable, for example, for pruning purposes.

In the following sections, the fundamental the RBF interpolation will be described. The finding of stationary points of surface using the RBF interpolation will be described in \cref{sec:ProposedApproach}. Moreover, the method, how the bindings between stationary points are searching, is introduced in this section. In the section \cref{sec:Experiments}, the results of our proposed algorithm will be presented. Finally, a final discussion of results will be performed.

\section{RBF Interpolation}\label{sec:RBFInterpolation}
In this section, the RBF interpolation method, recently introduced, e.g. in \cite{Skala2017}, \cite{Smolik2018}, and its properties are described.

We assume that we have an unordered dataset $\{\bm{x}_i\}_1^N\in\mathbb{E}^n$, where $n$ denotes the dimension of space and $N$ is the number of given points. Further, each point $\bm{x}_i$ from the dataset is associated with a vector $\bm{h}_i\in\mathbb{E}^p$ of the given values, where $p$ is the dimension of the vector, or a scalar value, i.e. $h_i\in\mathbb{E}^1$. In the following, we will deal with scalar data interpolation, i.e. the case when each point $\bm{x_i}$ is associated with a scalar value $h_i$ is considered. Our goal is determined the unknown function which is sampled at given points $\{\bm{x}_i\}_1^N$ by values $\{h_i\}_1^N$. For these purposes, it can be used the RBF interpolation which is based on the distance computation between two points $\bm{x}_i$ and $\bm{x}_j$ from the given dataset.

\begin{figure*}[htb]
\centering
\subfloat{\includegraphics[width = .32\linewidth]{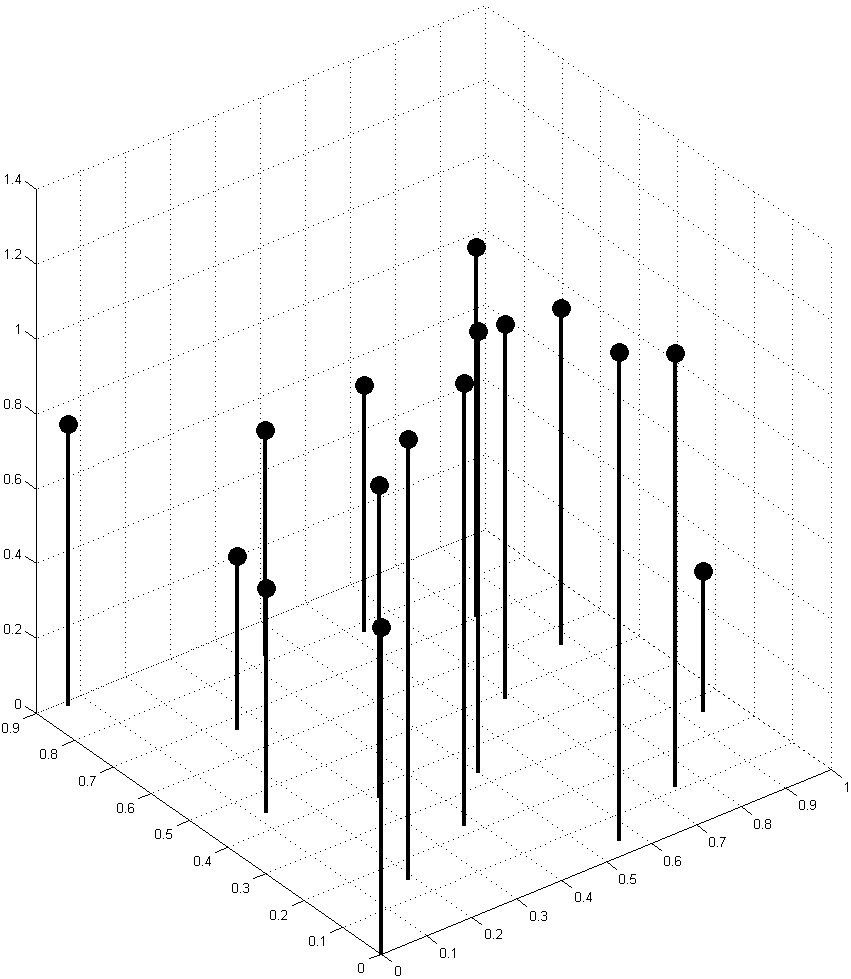}}\;
\subfloat{\includegraphics[width = .32\linewidth]{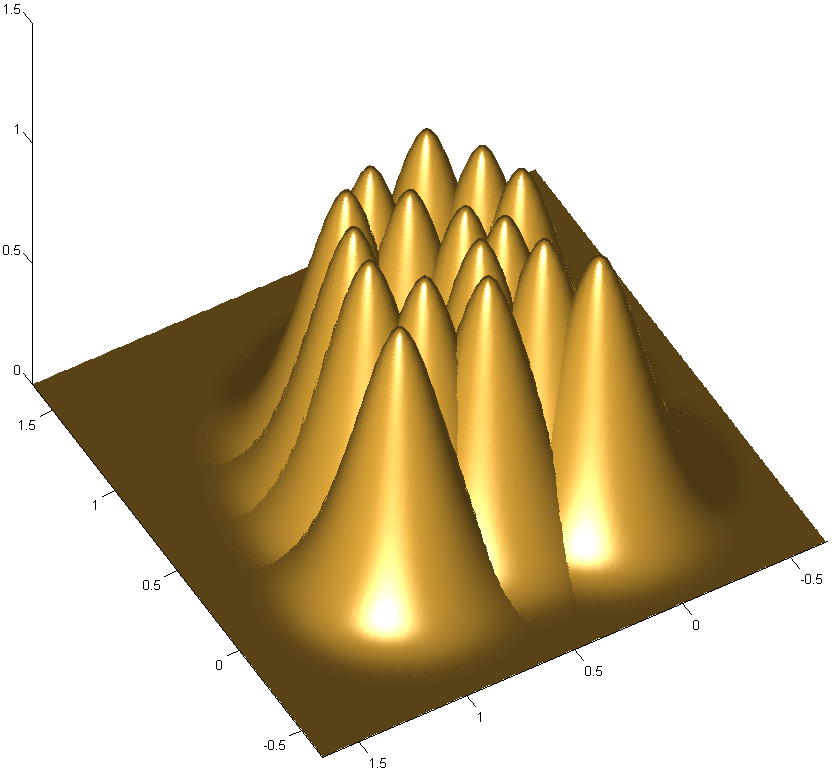}}\;
\subfloat{\includegraphics[width = .32\linewidth]{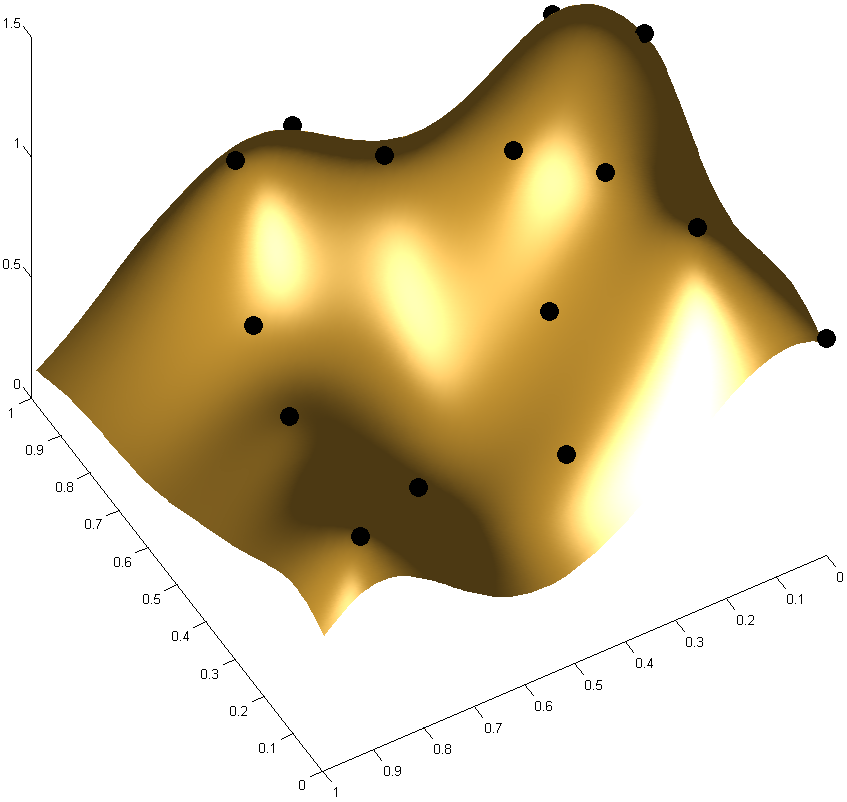}}\\
\caption{\label{fig:RBFinterp} Data values, the RBF collocation functions, the resulting interpolant.}
\end{figure*}

The interpolated value can be determined as:
\begin{equation}
\label{eq:interpFunc}
\displaystyle f(\bm{x}) = \sum_{j=1}^N c_j\phi(r_j)=\sum_{j=1}^N c_j\phi\left(\|\bm{x}-\bm{x}_j\|_2\right)\textrm{,}
\end{equation}
\noindent where the interpolating function $f(\bm{x})$ is represented as a sum of $N$ RBFs, each centered at a different data points $\bm{x}_j$ and weighted by an appropriate weight $c_j$ which has to be determined, see \cref{fig:RBFinterp}.

Applying (\ref{eq:interpFunc}) for all data points $\bm{x}_i, i = 1, \dots, N$, we get a linear system of equations:
\begin{equation}
\label{eq:interpLinEq}
\displaystyle h_i = f(\bm{x}_i) =\sum_{j=1}^N c_j\phi\left(\|\bm{x}_i-\bm{x}_j\|_2\right)\quad i = 1, \dots, N\textrm{.}
\end{equation}

The linear system of equations can be represented in a matrix form as:
\begin{equation}
\label{eq:matrixInterp}
\bm{Ac}=\bm{h}\textrm{,}
\end{equation}
\noindent where the matrix $\bm{A} = \left\{A_{ij}\right\} = \left\{\phi\left(\|\bm{x}_i-\bm{x}_j\|_2\right)\right\}$ is  $N\times N$ symmetric square interpolation matrix, the vector $\bm{c}=(c_1,\dots,c_N)^T$ is the vector of unknown weights and $\bm{h}=(h_1,\dots,h_N)^T$ is a vector of values in the given points. This linear system of equations can be solved by the Gauss elimination method, the LU decomposition, etc.

From the above, it can be seen that, in order to solve the interpolation problem, the distance matrix and a radial basis expansion are used.

\section{Proposed Approach}\label{sec:ProposedApproach}
In this section, determination of stationary points of the given dataset is described. Moreover, the approach includes the method for searching of bindings between stationary points because whole shape of stationary points may lie on the sampled surface.

\subsection{Piecewise approach for determination of stationary points}\label{sec:DetStatPts}
For simplicity we assume that we have given dataset $\{\bm{x}_i\}_1^N\in\mathbb{E}^2$ and each point $\bm{x}_i$ from this dataset is associated with a scalar value $h_i\in\mathbb{E}^1$. Further, for purposes of determination of stationary points, we assume that the given dataset contains the points on a $N_x\times N_y$ regular grid, where $\Delta x$ and $\Delta y$ are real numbers representing its grid spacing. Moreover, the row-major ordering of the given data is performed at first. After that, the piecewise approach is applied on the given data.

\begin{figure*}[htb]
\centering
\subfloat[\label{fig:piecewise}The grey area shows the all points from the given dataset which are interpolated by the RBF method during one step of piecewise approach. The hatched area illustrates the domain for which the stationary points of the given dataset are determined from the obtained RBF interpolation.]{\includegraphics[width = .48\linewidth]{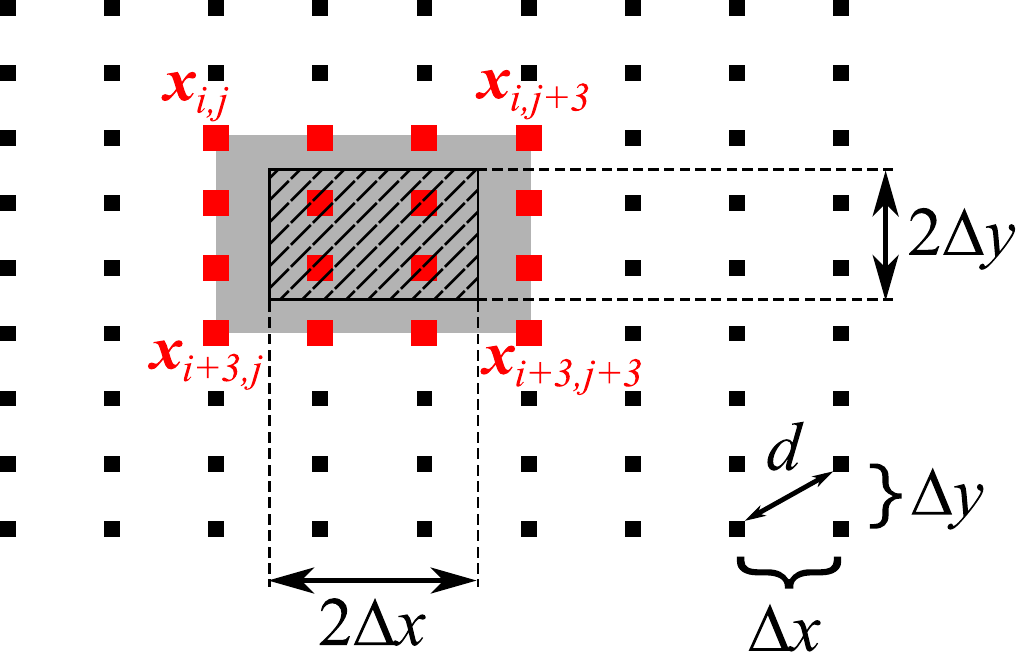}}\quad
\subfloat[\label{fig:piecewiseStep} Visualization of the stationary points reduction which is performed if the two points are identical or very close to identical. The green circle and red circle mark the stationary points which were determined from two different RBF interpolations and which were merged to one stationary point marked by yellow square.]{\includegraphics[width = .48\linewidth]{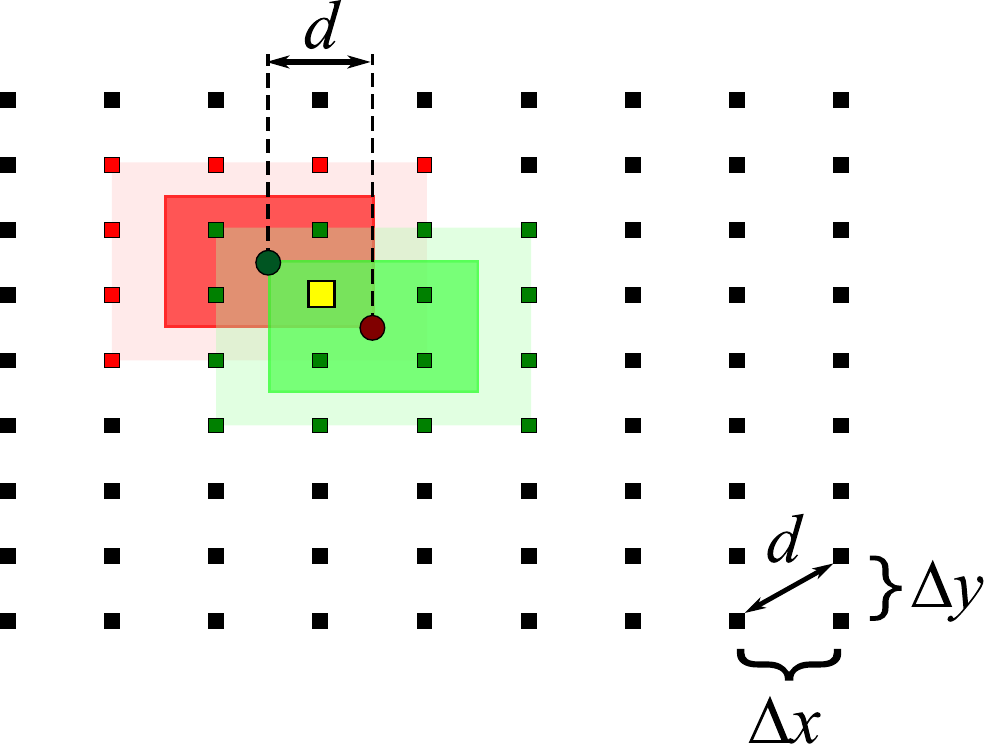}}\\
\caption{\label{fig:PWInterp} Proposed piecewise approach}
\end{figure*}

The process which is performed at each step of the piecewise approach is following. Every sixteen points $\{\bm{x}_m\}_1^{16}=\left\{\bm{x}_{i,j},\dots,\bm{x}_{i,j+3},\dots, \bm{x}_{i+3,j}, \dots,\bm{x}_{i+3,j+3} \right\}$ from the given dataset, where $i \in \{1,\dots, N_y-3\}$ denotes the row index and $j \in \{1,\dots, N_x-3\}$ denotes the column index, are interpolated by the RBF interpolation (\ref{eq:interpFunc}), i.e. the linear system (\ref{eq:matrixInterp}) has to be solved and the vector of weights $\bm{\hat{c}}=(c_1,\dots,c_{16})$ is computed. It mean that during one step of proposed approach, the RBF interpolation for $3\Delta x \times 3\Delta y$ area, where $\Delta x$ and $\Delta y$ are real numbers representing the input grid spacing, is performed, see \cref{fig:piecewise}. 

Then, the stationary points $\{\bm{s}_q\}$ of this interpolation function are determined using (\ref{eq:stacPts}). Specifically, for stationary points of the RBF interpolation function the nonlinear system of equations:
\begin{equation}
\label{eq:gradRBF}
\displaystyle \bm{0} = \sum_{m=1}^{16} c_m\frac{\phi'\left(\|\bm{x}-\bm{x}_m\|_2\right)}{\|\bm{x}-\bm{x}_m\|_2}\ast\left(\bm{x}-\bm{x}_m\right)\textrm{,}
\end{equation}
\noindent where $\phi'(r)$ is the derivation of RBF function $\phi$ with respect to variable $r$, $\ast$ denotes the element-wise multiplication and $\bm{\hat{c}}=(c_1,\dots,c_{16})$ is the vector of weights, has to be solved. The solution of (\ref{eq:gradRBF}), i.e. the stationary points $\{\bm{s}_q\}$ of the RBF interpolation, is searched for the domain defined as: 
\begin{equation}
\label{eq:AABB}
\begin{gathered}
\bm{x}_{i,j} + \bm{\varepsilon}_{min}\le\bm{s}_q\le\bm{x}_{i+3,j+3}- \bm{\varepsilon}_{max}\textrm{,}\\[0.2em]
\textstyle
\begingroup
 \bm{\varepsilon}_{min} = 
  \begin{cases} 
  \textstyle\left[\frac{\Delta x}{2},0\right]       & \textstyle\text{if } i=1\\[0.3em]
   \textstyle\left[0,\frac{\Delta y}{2}\right]       & \textstyle\text{if } j = 1\\[0.3em]
   \textstyle\left[\frac{\Delta x}{2},\frac{\Delta y}{2}\right] &  \textstyle\text{otherwise}
  \end{cases}\;
  \bm{\varepsilon}_{max} = 
  \begin{cases} 
   \textstyle\left[\frac{\Delta x}{2},0\right]       & \textstyle\text{if } i=N_y-3\\[0.3em]
   \textstyle\left[0,\frac{\Delta y}{2}\right]       & \textstyle\text{if } j=N_x-3\\[0.3em]
    \textstyle\left[\frac{\Delta x}{2},\frac{\Delta y}{2}\right] & \textstyle\text{otherwise}
  \end{cases}
  \endgroup
  \end{gathered}
\end{equation}
\noindent where $\Delta x$ and $\Delta y$ are real numbers representing the input grid spacing, $N_x$ indicates the number of grid column and $N_y$ is the number of grid rows, see \cref{fig:piecewise}, and the resulting set is added to the set of stationary points $\{\bm{s}_l\}$. It should be noted, that the values $\bm{\varepsilon}_{min}$ and $\bm{\varepsilon}_{max}$ include the correction for the boundary areas. 

The determination of stationary points of a function corresponds to the problem of finding critical points of the vector field, where the vector field is defined by eq. (\ref{eq:gradRBF}) for our purposes, and, therefore the method for determining critical points \cite{bhatia2014}, \cite{wang2018} may be used for obtaining the result.

The advantage of the above mentioned process is that the matrix $\bm{A}$ of the linear system (\ref{eq:matrixInterp}) for the RBF interpolation, is not dependent on the position of the given points (the matrix is dependent only on the distances between given points) and, therefore, this matrix is constant for all steps of piecewise approach. It should be noted that the approximation by a quadric surface could be used instead of the RBF interpolation, but the experimental results proved that this variant returns worse results in terms of stationary point locations.

The set of stationary points in the current form $\{\bm{s}_l\}$ may contain two identical points or points very close to identical. This problem is caused by the fact that one stationary point can be obtained from more RBF interpolations. The situation is illustrated in \cref{fig:piecewiseStep}. However, this problem can be solved by reduction of the set of stationary points. Then, the final set of stationary points $\left\{\bm{\sigma}_u\right\}$ of the given data is determined as follows. The subset $S_u$ of stationary points is removed from the unreduced set of stationary points $\{\bm{s}_l\}$. The points in the subset $S_u$ meet relation:
\begin{equation}
\label{eq:subsetStat}
S_u = \left\{\bm{s}_{k}:\|\bm{s}_{k}-\bm{s}_{1}\|\le d\right\}\textrm{,}
\end{equation}
\noindent where $d=\sqrt{(\Delta x)^2 +(\Delta y)^2}$ is the diagonal step in the regular grid, and the new stationary point is determined as a centroid of points from subset $S_u$:
\begin{equation}
\label{eq:newStat}
\bm{\sigma}_u =\frac{\sum\bm{s}_{k}}{|S_u|}\textrm{,}
\end{equation}
\noindent where $|S_u|$ is a number of points in the subset $S_u$. The process is repeated until the unreduced set of the stationary points is not empty.

The whole algorithm for determining the stationary points of the given dataset is summarized in \cref{alg:findStat}.

\renewcommand{\algorithmicrequire}{\textbf{Input:}}
\renewcommand{\algorithmicensure}{\textbf{Output:}}
\renewcommand{\algorithmicforall}{\textbf{for each}}   

\begin{algorithm*}[htb]
\floatname{algorithm}{Algorithm}
  \caption{Determination of the stationary points $\left\{\bm{\sigma}_u\right\}_1^{N_S}$.\label{alg:findStat}}
  \begin{algorithmic}[1]
  	\Require{given points $\{\bm{x}_i\}_1^N$ and their associated scalar values $\{h_i\}_1^N$, size of grid $N_x\times N_y$, grid spacing $\Delta x$ and $\Delta y$, used RBF $\phi$ and its shape parameter $\alpha$.}
	\Ensure{stationary points $\left\{\bm{\sigma}_u\right\}_1^{N_S}$}
	\State Row-major ordering the given points $\{\bm{x}_i\}_1^N$.
	\State $d=\sqrt{(\Delta x)^2 +(\Delta y)^2}$.
	\State Compute matrix $\bm{A}$ of linear system (\ref{eq:matrixInterp}) for the set of points $\left\{\bm{x}_{1},\dots,\bm{x}_{4},\bm{x}_{N_x+1},\dots,\bm{x}_{N_x+4}, \bm{x}_{2N_x+1},\dots,\bm{x}_{2N_x+4},\bm{x}_{3N_x+1},\dots,\bm{x}_{3N_x+4} \right\}$.
	\For{$i=1,\dots,N_y-3$}
		\For{$j=1,\dots,N_x-3$}
			\State $\bm{\hat{x}} = \left\{\bm{x}_{(i-1)N_x + j},\dots,\bm{x}_{(i-1)N_x + j+3},\dots, \bm{x}_{(i+2)N_x + j}, \dots,\bm{x}_{(i+2)N_x +j+3} \right\}$ 
			\State $\bm{\hat{h}} = \left\{h_{(i-1)N_x + j},\dots,h_{(i-1)N_x + j+3},\dots, h_{(i+2)N_x + j}, \dots,h_{(i+2)N_x +j+3} \right\}$ 
			\State Compute the vector of unknown weights $\bm{\hat{c}}$, eq. (\ref{eq:matrixInterp}), where $\bm{h}=\bm{\hat{h}}$.
			\State Compute the coefficients $\bm{\varepsilon}_{min}$ and $\bm{\varepsilon}_{max}$, eq. (\ref{eq:AABB}).
			\State Determine the stationary points $\{\bm{s}_{q}\}$ from eq. (\ref{eq:gradRBF}) in the domain (\ref{eq:AABB}).
			\State  $\{\bm{s}_l\} = \{\bm{s}_l\}\cup\{\bm{s}_{q}\}$
		\EndFor
	\EndFor
	\While{the set $\{\bm{s}_k\}$ is not empty}
		\State Find $S_u=\left\{\bm{s}_{k}: \|\bm{s}_{k}-\bm{s}_1\|_2\leq d\right\}$ in the set $\{\bm{s}_l\}$.
		\State Add the stationary point $\frac{\sum\bm{s}_{k}}{|S_u|}$ to the final set of stationary points $\left\{\bm{\sigma}_u\right\}$.
		\State Delete all points $\bm{s}_{k}\in S_u$ from the set $\{\bm{s}_l\}$.
	\EndWhile
  \end{algorithmic}
\end{algorithm*}

\subsection{Estimation of shape parameter for RBF interpolation}\label{sec:EstimationShape}
The piecewise RBF interpolation is used during the process of the determining the stationary points of the given dataset. Nevertheless, the quality of the resulting RBF interpolation strongly depends on the choice of the shape parameter $\alpha$. Therefore, in this section, the determination of suitable shape parameter $\alpha$ will be performed.

For the above mentioned process, the surface with the least possible tension is required, i.e. the surface must contain as little wavy as possible if the interpolated points allow it. It means that the shape parameter $\alpha$ has to be sufficiently large.

Therefore, for these purposes, we proposed and experimentally verified that shape parameter $\alpha$ is chosen so that the radius of circle of non-stationary inflection points of used RBF function $\phi(r)$ corresponds to the maximum distance of the interpolated points, within one step of proposed piecewise approach, which is $r=3d$, where $d =\sqrt{(\Delta x)^2 +(\Delta y)^2}$ is the diagonal step in the regular grid.

From this assumption, the following expression for shape parameter $\alpha$ was derived:
\begin{equation}
\label{eq:shape}
\alpha = \frac{\omega}{3d}\textrm{,}
\end{equation}
\noindent where $d =\sqrt{(\Delta x)^2 +(\Delta y)^2}$ is the diagonal step in the regular grid and $\omega$ is a constant parameter depending on the type of used RBF, see \cref{tab:Functions}.
\vspace{-1em}
\begin{table}[!h]
\centering
\caption{\label{tab:Functions} Different RBFs, their derivation $\phi'(r)$ and their parameter $\omega$, eq. (\ref{eq:shape}).}
\renewcommand{\arraystretch}{1.5}
\begin{tabular}{|c | c | c | c |}
\hline
{\textbf{RBF}} &{$\boldsymbol{\phi}(\bm{r})$}& {$\boldsymbol{\phi}\bm{'}(\bm{r})$} & {$\boldsymbol{\omega}$}\\ \hline \hline
{Gaussian RBF} & $\textstyle e^{-(\alpha r)^2}$& $\textstyle -2 \alpha^2 r e^{-(\alpha r)^2}$ & {$1/\sqrt{2}$} \\ \hline
{Inverse quadric} & $\textstyle\left(1+(\alpha r)^2\right)^{-1}$ &  $\textstyle -2 \alpha^2 r \left(1+(\alpha r)^2\right)^{-2}$ & {$1/\sqrt{3}$} \\ \hline
{Wendland's $\phi_{3,1}$} &  {$\textstyle(1-\mathnormal{\alpha} r)_{+}^4(4\alpha r + 1)$}& $\textstyle -20 \alpha^2 r(1-\mathnormal{\alpha} r)_{+}^3$ & {$1/4$} \\ \hline
\end{tabular}
\end{table}
\vspace{-2em}
\subsection{Searching of bindings between stationary points}\label{sec:Binding}
It is possible that the given surface does not contain only isolated stationary points, but the curves of stationary points, such as line segments, circles, parabolas or some other shapes, can lie on the given surface. Therefore, the method for searching of bindings between stationary points will be described.

\begin{figure}[t]
\centering\includegraphics[width = 0.48\linewidth]{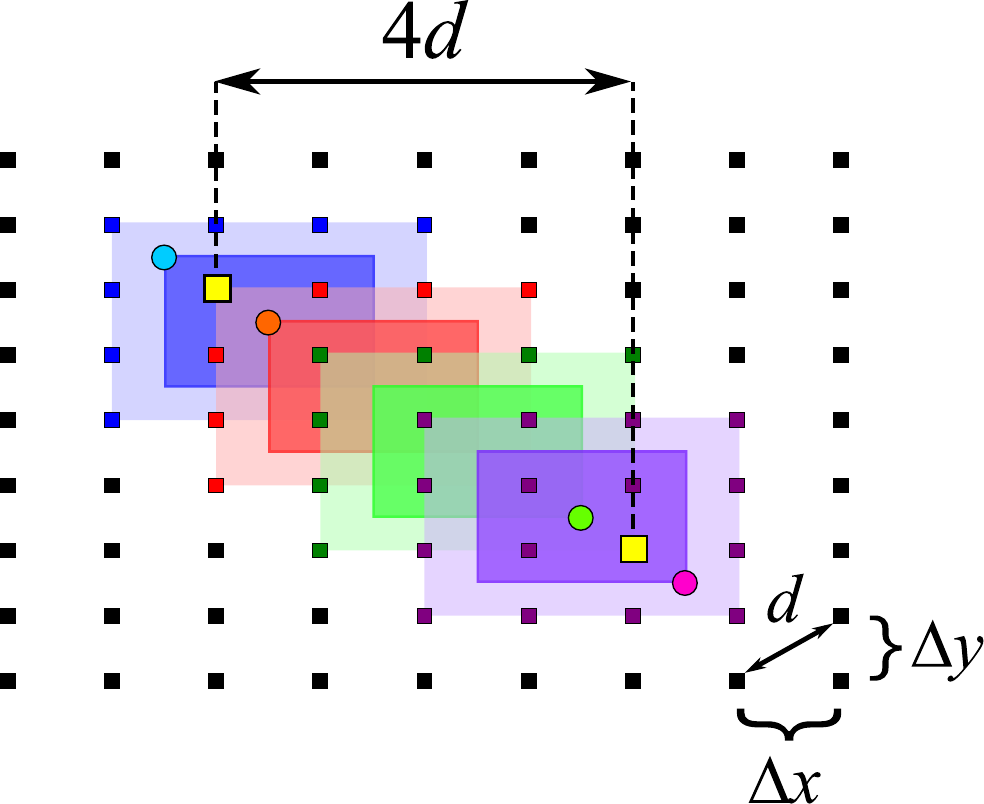}
\caption{\label{fig:piecewiceFace} The figure shows the worst case in which two stationary points (yellow squares) still lie on the same curve of stationary points, i.e. the distance between two stationary points is maximal possible distance. Moreover, the reduction of stationary points is again shown.}
\end{figure}

At the beginning, the maximal possible distance $\delta_{max}$ of two stationary points for which these stationary points still lie on the same curve has to be established. The situation of the worst case is illustrated in \cref{fig:piecewiceFace}. 
In this figure, it can be seen four subdomains of the piecewise approach and for each of them, the one stationary point is indicated using circle mark. Based on eq. (\ref{eq:subsetStat}), the stationary points of blue and red subdomains are reduced and are replaced by their centroid. The same case occurs for the green and purple subdomain. New stationary points which are obtained after the reduction are indicates by yellow squares in the figure. It is also obvious that the distance of these two stationary points, which is also the maximum possible distance $\delta_{max}$, is:
$$\delta_{max}=4d\textrm{,}$$
\noindent where $d =\sqrt{(\Delta x)^2 +(\Delta y)^2}$ is the diagonal step in the regular grid.

Now, the stationary points $\left\{\bm{\sigma}_u\right\}$ of the given dataset are sequentially processed by following. For the current stationary point $\bm{\sigma}_u$, the all stationary points $\left\{\bm{\sigma}_w\right\}$ which lying in the distance $\delta_{max}$ are determined:
\begin{equation}
\label{eq:bind}
 \left\{\bm{\sigma}_w\right\} = \left\{\bm{\sigma}_w:\|\bm{\sigma}_w-\bm{\sigma}_u\|\le \delta_{max}\right\}\textrm{.}
\end{equation}
\noindent If no stationary point is found, then the stationary point $\bm{\sigma}_u$ is isolated. Otherwise, the binding $f_v= \left\{\bm{\sigma}_w\right\}\cup\left\{\bm{\sigma}_u\right\}$ is obtained and newly added stationary points are processed in the same way. Finally, the result of this approach is the set of points described the curve of stationary points. This procedure is repeated until the all stationary points $\left\{\bm{\sigma}_u\right\}$ are processed.

One of the possible solution of this problem is the $kd-$tree which can be simply applied for purposes of searching of bindings between stationary points.

\section{Experimental Results}\label{sec:Experiments}
In this section, the experimental results for our proposed approach will be presented and their comparison with the exact stationary points which were determined analytically from the sampling function will be made. The implementation was performed in Matlab. In addition, different radial basis functions have been used, see \cref{tab:Functions}.

For purposes of our experiments, a uniform distribution of points on a rectangular domain was used for the testing data. Thus, the given dataset contains $120\times120$ points uniformly distributed in the interval $[x_{min},x_{max}]\times[y_{min}, y_{max}]$, where the values $x_{min}$, $x_{max}$, $y_{min}$ and $y_{max}$ are chosen based on the used sampling function, see (\ref{eq:TestF1})~-~(\ref{eq:TestF8}) and (\ref{eq:TestFL})~-~(\ref{eq:TestFT}). Moreover, each point from this dataset is associated with a function value of the selected sampling function at this point.

\subsection{Comparison of determined stationary points with exact stationary points}\label{sec:Comparison}
In this section, the results for datasets whose stationary points do not contain mutual bindings, i.e. all stationary points are isolated, will be presented. The sampling functions $f_1$ (\ref{eq:TestF1}) and $f_2$ (\ref{eq:TestF8}), which were defined in \cite{franke1979}, fulfill these properties.

\begin{subequations}
\label{eq:exp1}
\begin{align}
\begin{split}
\label{eq:TestF1} f_1(\bm{x}) & = \frac{3}{4}e^{-\frac{(9x_1-2)^2}{4}-\frac{(9x_2-2)^2}{4}} + \frac{3}{4} e^{-\frac{(9x_1+1)^2}{49}-\frac{(9x_2+1)}{10}} \\
&\;+ \frac{1}{2} e^{-\frac{(9x_1-7)^2}{4}-\frac{(9x_2-3)^2}{4}} - \frac{1}{5} e^{-(9x_1-4)^2-(9x_2-7)^2}
\end{split}\quad \bm{x}\in[0,1]\times[0,1]\\
\label{eq:TestF8} f_2(\bm{x}) & =\sin{\left(3\cdot x_1\right)}\cdot\cos{\left(3\cdot x_2\right)}\qquad\qquad\qquad\qquad\quad \bm{x}\in[-2,2]\times[-2,2]
\end{align}
\end{subequations}

\begin{figure*}[htb]
\centering
\subfloat[\label{fig:F1Gauss} Dataset sampled from $f_{1}$, eq.(\ref{eq:TestF1})]{\includegraphics[width = .48\linewidth]{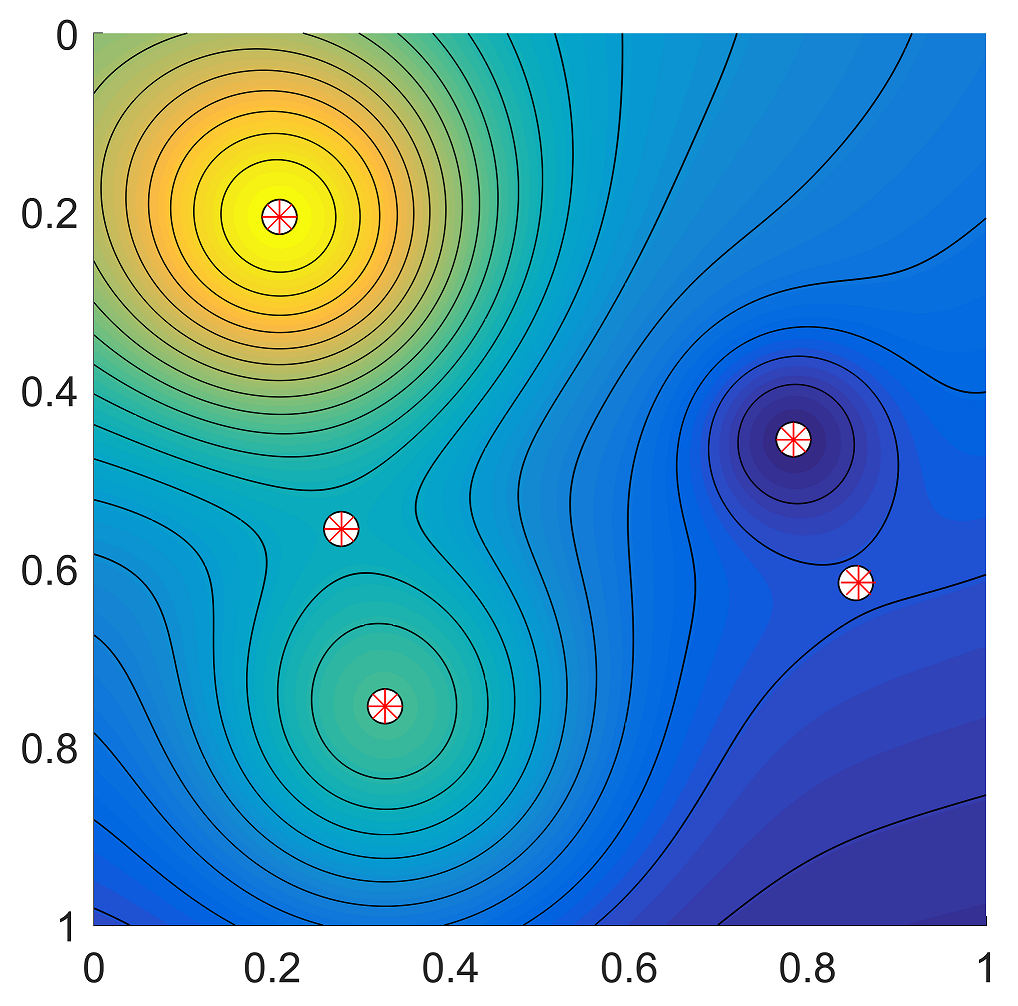}}\;
\subfloat[\label{fig:F8Gauss} Dataset sampled from $f_{2}$, eq.(\ref{eq:TestF8})]{\includegraphics[width = .48\linewidth]{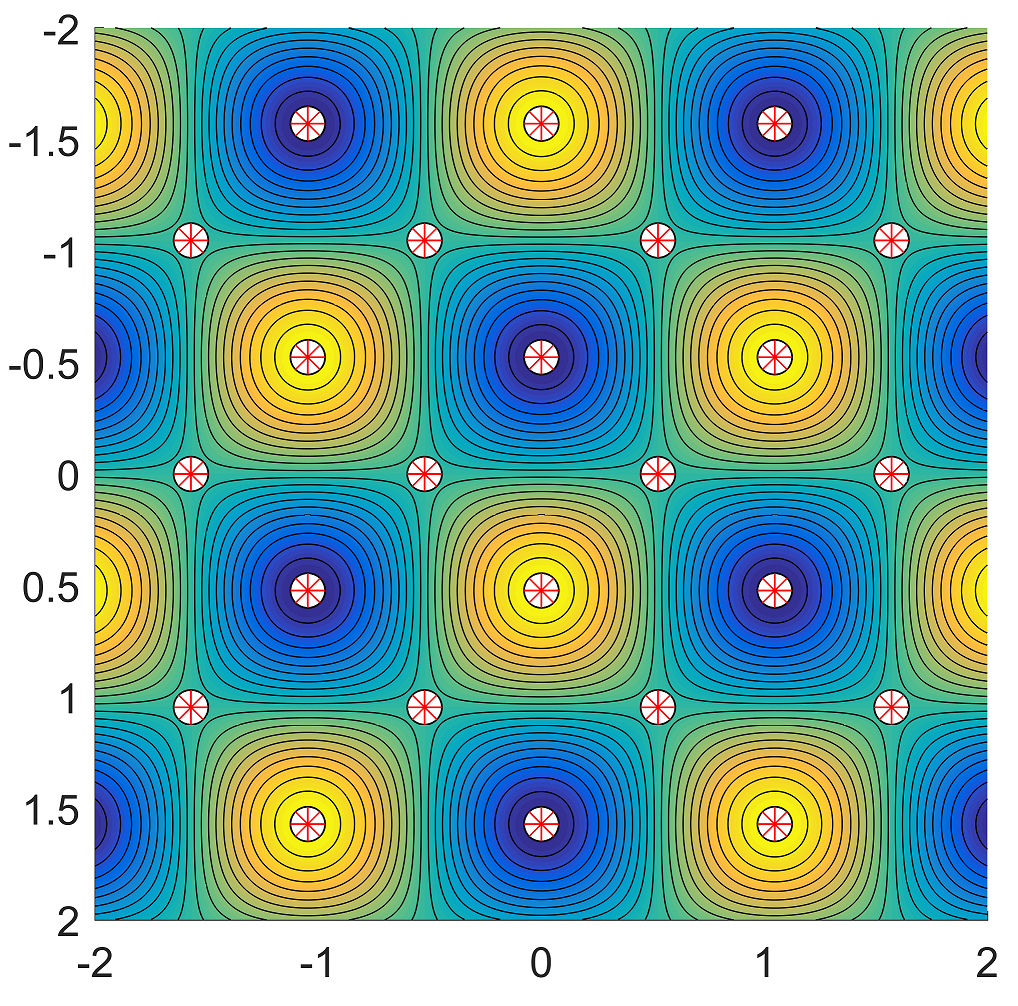}}\\
\caption{\label{fig:exp1} The white circles indicate the stationary points of the given dataset that are obtained using the proposed approach when the RBF interpolation used the Gaussian RBF. The tested dataset contains $120\times 120$ points. The red asterisks ($*$) denote the exact positions of the stationary points of the appropriate function. Furthermore, the contour map of given dataset is shown.}
\end{figure*}

\Cref{fig:F1Gauss} presents the results for the dataset in which each point is associated with a value from the $f_1$ function (\ref{eq:TestF1}) when the Gaussian RBF has been used for the piecewise RBF interpolation. Using our proposed approach, five isolated stationary points which are marked by white circles were found for this dataset. The exact stationary points of $f_1$ function (\ref{eq:TestF1}) are shown using the red asterisks ($*$).

The results for the dataset in which each point is associated with a value from the $f_2$ function (\ref{eq:TestF8}), when the Gaussian RBF has been used for the piecewise RBF interpolation, are presented in \cref{fig:F8Gauss}. Twenty four isolated stationary points which are represented by white circles were found for this dataset using our proposed approach. The exact stationary points of $f_2$ function (\ref{eq:TestF8}) are again shown using the red asterisks ($*$).

It can be seen that obtained results for both datasets correspond to the stationary points calculated analytically from the sampling functions. Moreover, it should be noted, that the same results were obtained even when other RBF function, see \cref{tab:Functions}, were used for the piecewise RBF interpolation.

\subsection{Bindings between stationary points}\label{sec:ExperimentBinding}
In this section, the results for datasets whose stationary points contains mutual bindings will be presented.  The four following sampling functions  (\ref{eq:TestFL})~-~(\ref{eq:TestFT}) fulfill these properties.

\begin{subequations}
\label{eq:exp2}
\begin{align}
\label{eq:TestFL} f_{11}(\bm{x}) & = -\left(x_1-x_2\right)^2\textrm{,}&\quad \bm{x}\in[-1,1]\times[-1,1]\\ 
\label{eq:TestFP} f_{12}(\bm{x}) & =\sin{\left(x_1+x_2^2\right)}\textrm{,}&\quad \bm{x}\in[-3,3]\times[-2,2]\\
\label{eq:TestFO} f_{13}(\bm{x}) & =\sin{\left(3\pi\left(\sqrt{x_1^2+x_2^2}+0.25\right)\right)}\textrm{,}&\quad \bm{x}\in[-1,1]\times[-1,1]\\
\label{eq:TestFT} f_{14}(\bm{x}) & =-2\cdot\left(x_1^2-x_2^2\right)^2+1\textrm{,}&\quad \bm{x}\in[-1,1]\times[-1,1]
\end{align}
\end{subequations}

At the beginning, it should be noted that the white solid line indicates the curve of stationary points obtained for the given dataset using our proposed approach and the isolated stationary point obtained using our proposed approach is marked by the white circle. The red dashed line indicates the curve of stationary points calculated analytically from the given sampling function and the isolated stationary point calculated analytically from the given sampling function is represented by the red asterisk ($*$).

\begin{figure}[!b]
\centering
\subfloat[\label{fig:FLGauss} Dataset sampled from $f_{11}$, eq.(\ref{eq:TestFL})]{\includegraphics[width = .48\linewidth]{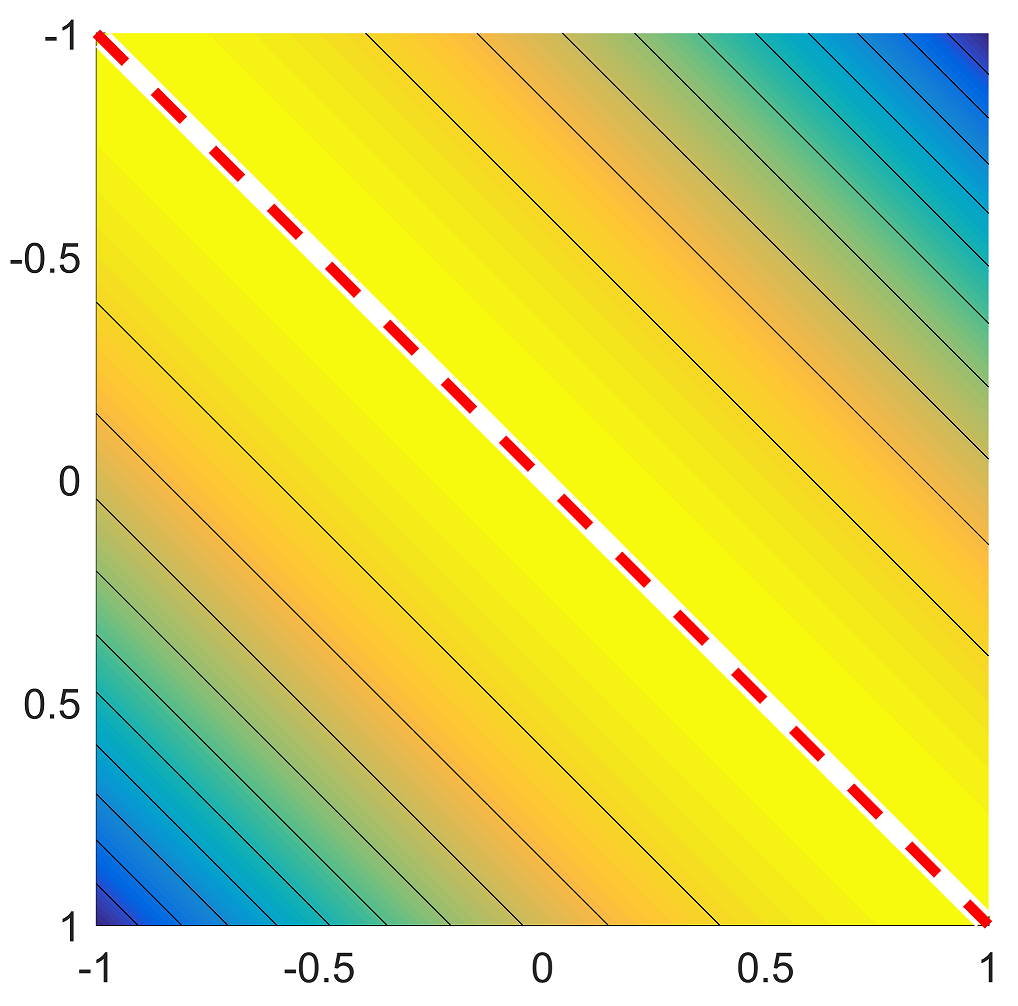}}\;
\subfloat[\label{fig:FPGauss} Dataset sampled from $f_{12}$, eq.(\ref{eq:TestFP})]{\includegraphics[width = .48\linewidth]{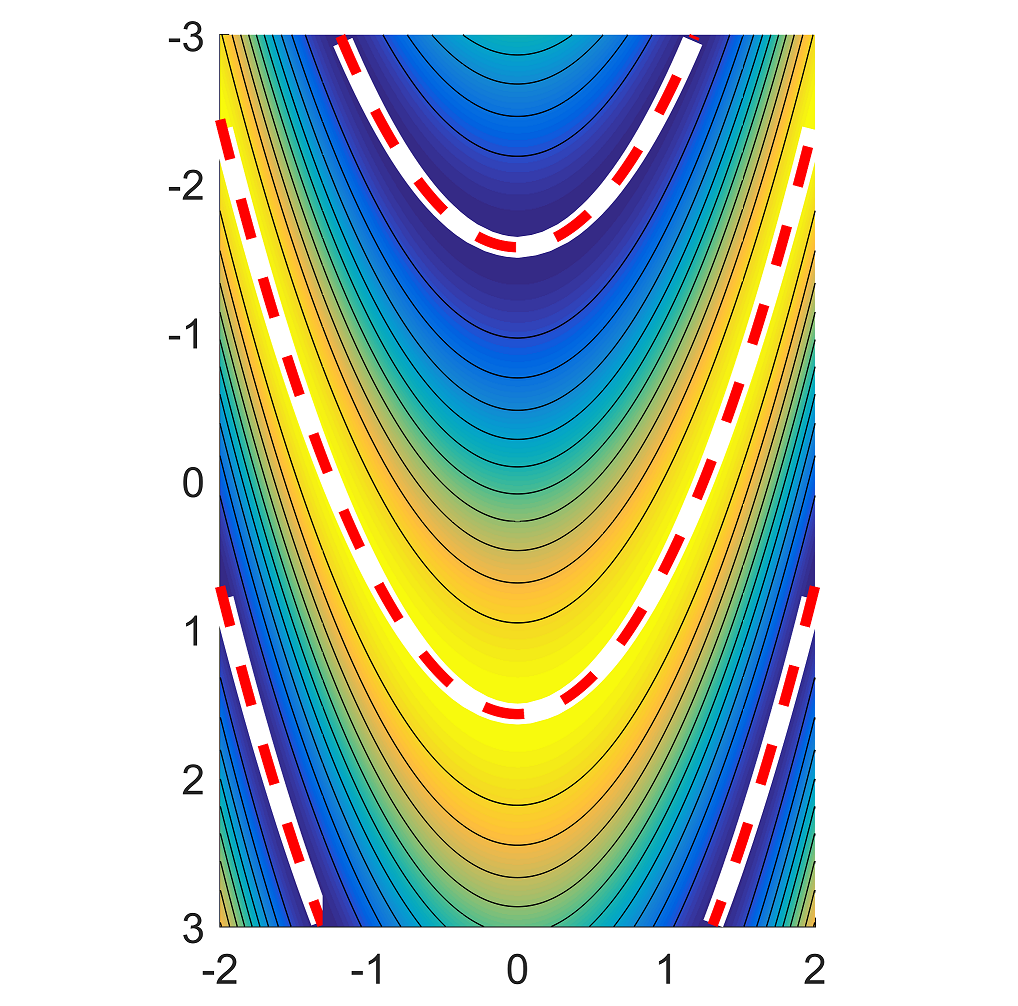}}\\
\subfloat[\label{fig:FOGauss} Dataset sampled from $f_{13}$, eq.(\ref{eq:TestFO})]{\includegraphics[width = .48\linewidth]{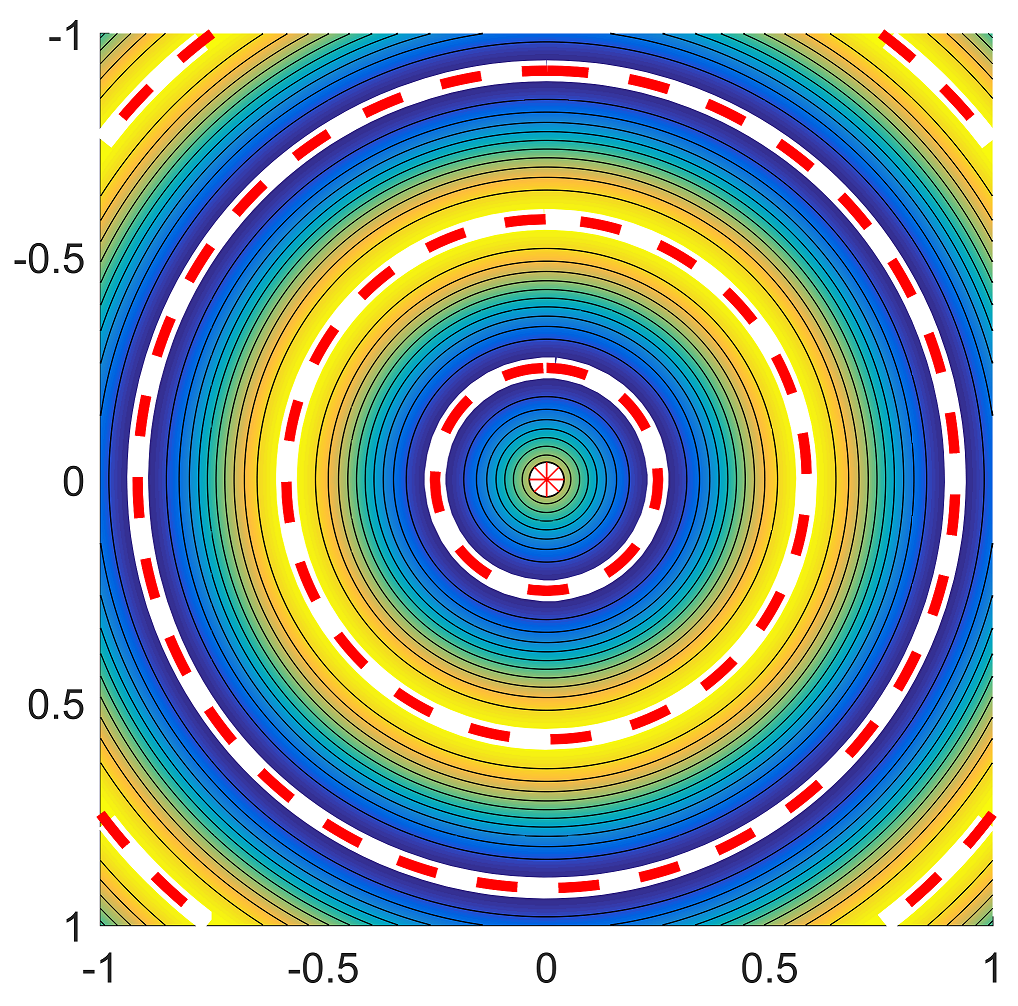}}\;
\subfloat[\label{fig:FTGauss} Dataset sampled from $f_{14}$, eq.(\ref{eq:TestFT})]{\includegraphics[width = .48\linewidth]{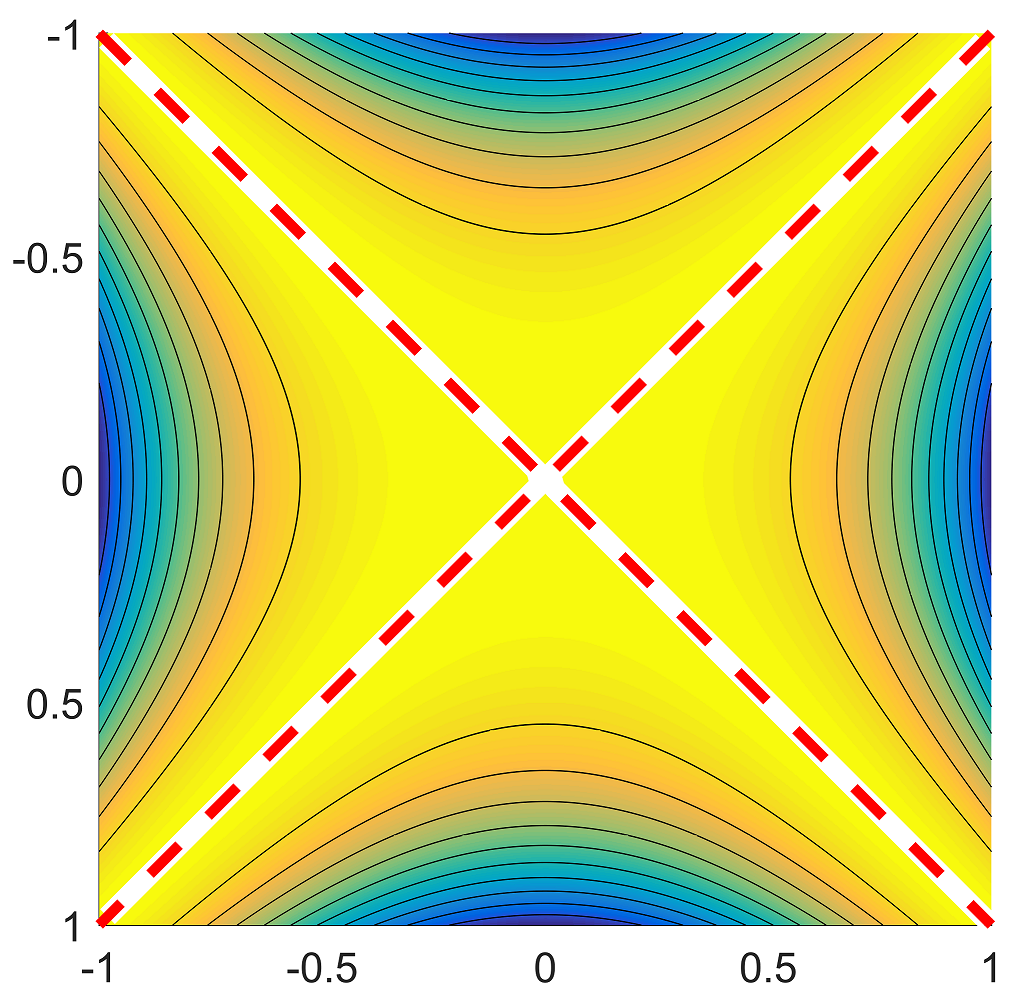}}\\
\caption{\label{fig:exp2} The white solid lines indicate the curves of stationary points of the given dataset that are obtained using the proposed approach when the RBF interpolation used the Gaussian RBF. The tested dataset contains $120\times 120$ points. The red dashed lines denote the exact curves of the stationary points of the appropriate function. Furthermore, the contour map of given dataset is shown.}
\end{figure}

\Cref{fig:FLGauss} presents the results for the dataset in which each point is associated with a value from the $f_{11}$ function (\ref{eq:TestFL}) when the Gaussian RBF has been used for the piecewise RBF interpolation. Using our proposed approach, one curve of stationary points, specifically the line segment, were found for this dataset. This result coincides with the result obtained using analytically approach.

The results for the dataset in which each point is associated with a value from the $f_{12}$ function (\ref{eq:TestFP}), when the Gaussian RBF has been used for the piecewise RBF interpolation, are presented in \cref{fig:FPGauss}. For this dataset, the four curves of stationary points, specifically two parabolas and two segments of parabola, were found using our proposed approach.

\Cref{fig:FOGauss} presents the results for the dataset in which each point is associated with a value from the $f_{13}$ function (\ref{eq:TestFO}) when the Gaussian RBF has been used for the piecewise RBF interpolation. Using our proposed approach, seven curves of stationary points, specifically three circles and four arcs,  and one isolated stationary point were found for this dataset.

The results for the dataset in which each point is associated with a value from the $f_{14}$ function (\ref{eq:TestFT}), when the Gaussian RBF has been used for the piecewise RBF interpolation, are presented in \cref{fig:FTGauss}. For this dataset, the two curves of stationary points, specifically two line segments, were found using our proposed approach.

For all mentioned experiments, it can be seen that the results obtained using our proposed approach correspond to results obtained using analytically approach. Moreover, it should be again noted that the same results were obtained even when other RBF function, see \cref{tab:Functions}, were used for the piecewise RBF interpolation.

\section{Conclusion}\label{sec:Conclusion}
In this paper, a new approach for determination of stationary points of given sampled surface without knowledge of the sampling function is presented. The proposed method is based on the piecewise RBF interpolation of the given dataset. Moreover, the proposed approach includes the method of detecting the bindings between the found stationary points, i.e. the approach is able to associate the points from the same curve of stationary points.

The experiments proved that the stationary points determined by our proposed approach coincide with the exact stationary points which were determined analytically form the sampling function.

The results of the proposed approach can, for example, be used for determination of the set of reference points for the RBF approximation which enable appropriate compression of given dataset. The knowledge of the bindings between stationary points is possible to use for pruning the subset of related stationary points to the required number of points on the appropriate curve of stationary points. 

In the future work, the proposed approach can be generalized for scattered data using the $k$-nearest neighbors algorithm.
\subsubsection*{Acknowledgments.}
The authors would like to thank their colleagues at the University of West Bohemia, Plze\v{n}, for their discussions and suggestions, and the anonymous reviewers for their valuable comments. Special thanks belong to Jan Dvorak, Lukas Hruda and Martin \v{C}ervenka for their independent experiments and valuable comments. The research was supported by the Czech Science Foundation GA\v{C}R project GA17-05534S and partially supported by the SGS~2016-013 project.

\bibliographystyle{splncs}
\bibliography{paper}

\end{document}